\documentclass{atmp}

\usepackage{cite}
\usepackage[all]{xy}
\usepackage{color}

\usepackage{epsfig}

\copyrightnotice{2006}{10}{1}{12}

\def\hpic #1 #2 {\mbox{$\begin{array}[c]{l} \epsfig{file=#1,height=#2} \end{arr\
ay}$}}
\def\vpic #1 #2 {\mbox{$\begin{array}[c]{l} \epsfig{file=#1,width=#2} \end{arra\
y}$}}

\newcommand{\aco}[2]{A^{[#1]}_{#2}}
\newcommand{\yco}[2]{Y^{[#1]}_{(#2)}}
\def\HOM{\mathop{{\rm Hom}}\nolimits}

\def\Z{\mathbb Z}
\def\C{\mathbb C}
\def\R{\mathbb R}
\def\P{\mathbb P}
\def\Q{\mathbb Q}

\def\O{\mathcal O}
\def\L{\mathcal L_K}

\def\RR{\mathcal R}

\numberwithin{equation}{section}

\begin{document}

\title[Seidel's mirror map for abelian varieties]{Seidel's
mirror map for abelian varieties}

\arxurl{math.SG/0512229a}

\author[Marco Aldi and Eric Zaslow]{Marco Aldi and Eric Zaslow}

\address{Department of Mathematics, Northwestern University,\\ 2033, Sheridan Road, Evanston,
IL 60208, USA}



\begin{abstract}
We compute Seidel's mirror map for abelian varieties by
constructing the homogeneous coordinate rings from the Fukaya
category of the symplectic mirrors. The computations are
feasible, as only linear holomorphic disks contribute to the
Fukaya composition in the case of the planar Lagrangians used.
The map depends on a symplectomorphism $\rho$ representing the
large complex structure monodromy.  For the example of the
two-torus, different families of elliptic curves are obtained by
using different $\rho$'s which are linear in the universal cover.
In the case where $\rho$ is merely affine linear in the universal
cover, the commutative elliptic curve mirror is embedded in
noncommutative projective space. The case of Kummer surfaces is
also considered.
\end{abstract}

\maketitle

\section{Introduction}\label{intro}

In \cite{Z}, we constructed Seidel's mirror map for the two-torus.
Starting from the Fukaya category of a symplectic two-torus $X$,
we computed the homogeneous coordinate ring ${\mathcal R}$ of the
mirror elliptic curve $Y.$ The proof of Kontsevich's conjecture
for the elliptic curve \cite{Pol-Z} allowed us to do so.  Namely,
$D{\rm Fuk}(X) \cong D(Y)$ implies that ${\mathcal R}$ is
computable on $X$ alone, i.e.,\vadjust{\pagebreak} we have
\begin{equation}\label{maineq}
\mathcal R = \bigoplus_{k=0}^\infty \Gamma(\mathcal O_Y(k)) =
\bigoplus_{k=0}^\infty \HOM_{D{\rm Fuk}(X)}(\psi(\mathcal O),
\psi(\mathcal O(k))),
\end{equation}
where $\psi$ is the equivalence of categories (see also
\cite{AB}).  In fact, in the above, $\mathcal O$ can be replaced
by any line bundle $\mathcal L.$ As the mirror of the
automorphism --- $\otimes {\mathcal O}_Y(1)$ is known to be the
symplectomorphism $\rho$ effecting the monodromy around the large
complex structure limit, the right-hand side can be computed
entirely from the Fukaya category, once we choose any Fukaya
object mirror to some line bundle. In this paper, we extend the
computation of \cite{Z} to higher-dimensional abelian varieties.
These cases were treated in \cite{fukabel} and \cite{kontsoib},
where partial results were obtained toward Kontsevich's
conjecture.  Those results imply the existence of some of the
findings in this paper, although our methods are more direct. We
find that the homogeneous coordinate ring of the mirror abelian
variety is described by the computation on the right-hand side of
equation~(\ref{maineq}). Specifically, we find that the theta
relations of the mirror abelian variety are obeyed by the
intersection points mirror to the theta functions. Computations
are made feasible by the fact that the objects in the Fukaya
category are all linear planes in the universal cover.

In addition, we explore the dependence of the homogeneous
coordinate ring on $\rho.$  Specifically, different choices for
$\rho$ yield different families of abelian varieties. When
$\rho(x,y)$ is not strictly linear but contains a translation, the
mirror elliptic curve sits as a commutative variety inside a
noncommutative projective space, not unlike the situation found in
\cite{AKO}. The Hesse family of elliptic curves was found to be
dual to $\rho(x,y) = (x,y+3x)$, in \cite{Z}. If~we choose
$\rho(x,y) = (x,y+x)$, we find the universal family of elliptic
curves inside weighted projective space ${\mathbb P}(1,2,3)$
corresponding to a quasihomogeneous coordinate ring with
corresponding weights. In the four-dimensional case, $\Z/2\Z$
invariant maps on the four torus provide a mirror map for
\hbox{Kummer} K3 surfaces.

\section{Definitions}\label{sec1}
Our main objects of study are symplectic tori with B-field
admitting the structure of Lagrangian brane fibration with
section (see \cite{BOOK,Pol}). This motivates the following
definitions. Let $V$ be a real vector space of dimension $2n$
together with a complexified symplectic form $i\omega+ b$ and a
basis $l_1 \ldots , l_n, l'_1, \ldots, l'_n$ for which
$(i\omega+b)(l'_i,l_j)=iM_{ij}+B_{ij}$ or, in matrix form,
\[
\omega = \left(
\begin{matrix}
0  & -M^{{\rm T}} \\
M & 0
\end{matrix}
\right)
\quad
b =
\left(
\begin{matrix}
0 & -B^{{\rm T}}\\
B & 0
\end{matrix}
\right)
\]
where $M = (M_{ij})$ and $B=(B_{ij})$ are invertible $n\times n$
matrices. Moreover,~let
\[
\Lambda_0:= \Z \langle l_1,\ldots, l_n \rangle  , \quad \Lambda'
:= \Z \langle l'_1, \ldots, l'_n \rangle .
\]
Consider the Lagrangian subspaces $L_0 := \Lambda_0 \otimes \R$,
$ L':= \Lambda' \otimes \R$. Suppose that there is a lattice
homomorphism $f: \Lambda_0 \to \Lambda'$ such that its
$\R$-linear extension (also denoted by $f$) satisfies $f(l_i) =
\sum_j N_{ji} l'_j$ and let $N := (N_{ij})$. Then the following
are equivalent
\begin{enumerate}
\item[(1)]
The graph $L_f := \{ l + f(l) | l \in L_0 \}$ is Lagrangian

\item[(2)]
$N^{\rm T} M$ and $N^{\rm T} B$ are symmetric

\item[(3)]
$i \omega + b$ admits the compatible complex structure
\[
J_f := \left(
\begin{matrix}
0 & - N^{-1} \\
N & 0
\end{matrix}
\right)
\]

\item[(4)]
$i \omega + b$ admits the symplectomorphism
\[
\rho_f := \left(
\begin{matrix}
{\rm Id} & 0 \\
N & {\rm Id}
\end{matrix}
\right).
\]
\end{enumerate}
Notice that if $f$ satisfies these equivalent conditions, then so
does $kf$ for all $k\neq 0$ and  $\rho_{(kf)} = \rho_f^k$.  In
more general settings, we drop the subscript of $\rho.$

Consider now the quotient torus $T := V / (\Lambda_0 \oplus
\Lambda')$. Each $L_k := L_{kf} = \rho_f^k L_0 $ has rational
slope with respect to the lattice and so it descends to a
Lagrangian subtorus, also denoted by $L_k$. We are interested in
the (derived) Fukaya subcategory generated by these Lagrangian
subtori. Define
\[
\Lambda(L_k)_0 := f^{-1}\!\left(\frac{1}{k} \Lambda'\right),\quad
K(L_k)_0 := \frac{\Lambda(L_k)_0}{\Lambda_0}.
\]
Then for $k_1\neq k_2$, we have the following isomorphisms of vector spaces
\begin{align*}
\HOM (L_{k_1},L_{k_2}) & \cong  \langle L_{k_1} \cap L_{k_2}
\rangle \cong \langle \{ (l,l')  | l' = k_1f (l) = k_2f(l) \,{\rm
mod}\,  \Lambda'\} \rangle
\\
& \cong  \langle K(L_{k_2-k_1})_0 \rangle \cong \langle L_0 \cap
L_{k_2-k_1} \rangle \cong \HOM  (L_0, L_{k_2-k_1}).
\end{align*}
Since $\rho_f^{k_2-k_1}$ is a symplectomorphism, the first term and last one are isomorphic as whatever structure the Fukaya category is enriched over. Using these isomorphisms, the composition
\begin{equation} \label{m}
m{:}\ \HOM (L_0, L_{k_1}) \otimes \HOM  (L_0, L_{k_2})
\longrightarrow \HOM (L_0, L_{k_1+k_2})
\end{equation}
is well defined.

\section{Only planar disks contribute}

The results in this section are not new. See for example
\cite{fukabel}.

All compositions relevant to our computation involve calculating
holomorphic maps from a disk with three marked points on the
boundary. The~three intervals on the boundary between pairs of
marked points must be~mapped to three Lagrangians, with the points
sent to intersections of those Lagrangians. In our case, the three
Lagrangians are all related to the base of the Lagrangian
fibration $L_0$ by the symplectomorphism $\rho_f$ linear in the
periodic coordinates, thus they all lift to planes in the
universal cover~$V$.

Since the image of a holomorphic map from a disk is connected, in
the universal cover, only one preimage of a Lagrangian will be
relevant to any given disk. As a result, we may perform our
calculation in the universal cover. Further, we may define
coordinates such that the preimages of $L_0$ and $L_{k_1}$
intersect at the origin $Y_0 = 0\in V$. Let $Y_1$  be the second
vertex, where $Y_1\in L_0 \cap L_{k_2-k_1}$. Then the minimal
choice for the third vertex is $Y_2 = ((k_2-k_1)/k_2) Y_1 \in L_0
\cap L_{k_2}$. With respect to the complex structure $J_f$, the
holomorphic map $g{:}\ D\to V$ may be constructed as
\[
 z_m (y)= a_m y \quad m\ge 1
\]
where $a_m$ are the coefficients in the decomposition $Y_1 =
\sum_{m=1}^n a_ m l_m$, and $y= y_1 + i y_2$ is the holomorphic
coordinate on a triangular domain (disk) $D$ in the complex plane
defined by the vertices $0$, $1$ and $({k_2}/({k_2 - k_1}))(1 +
k_1 i)$. The three legs of this triangle have slope $0$, $k_1$
and $k_2$, and we call these boundary intervals $C_0$, $C_1$ and
$C_2$, respectively. The map $g$ is easily shown to satisfy all
boundary conditions. Since the universal cover $V$ is a vector
space, any other holomorphic map $g'$ satisfying the boundary
conditions  may be expressed as $g' = g + \xi$, with $\xi$ a
vector-valued holomorphic function. Further, since $g$ and $g'$
satisfy the same affine linear boundary conditions, the boundary
conditions on $\xi$ are strictly linear. Explicitly, we have
\[
C_0{:}\ \overline \xi_j = \xi_j ; \qquad C_1{:}\ \overline \xi_j =
\mu_{\alpha(k_1)} \xi_j; \quad \overline \xi_j = \mu_{\alpha(k_2)}
\xi_j, \enspace j = 1,2,
\]
where we have defined $\alpha(n) = 1 + ni$ and $\mu_\alpha =
\overline \alpha / \alpha$. Note that $|\mu_\alpha| =1$. Consider
$\xi \equiv \xi_j$, where $j$ could be $1$ or $2$. Since $\xi$ is
homolorphic, it has a convergent power series expansion $\xi =
\sum_n c_n y^n$, and the boundary condition on $C_0$ requires
that the $c_n$ are real. For the boundary $C_1$, put $\alpha
\equiv \alpha(k_1) = 1 +  k_1 i$ and note that on $C_1$ we may
write $y= s \alpha$, where $s$ is real and runs from $0$ to
$k_2/(k_2-k_1)$. Then we write the boundary conditions as
\[
\overline \xi \mid_{C_1} =  \sum_{n=0}^\infty c_n (s\overline
\alpha)^n  = \mu_\alpha \xi\mid_{C_d} = \left(\frac{\overline
\alpha}{\alpha}\right) \sum_{n=0}^\infty c_n (s\alpha)^n .
\]
Equating powers of $s$ requires that for all $n$ such that $C_n$
is nonzero, we have $(\mu_\alpha)^{n-1} = 1$. If $n \neq 1$ we
learn that $\mu_\alpha$ must be a root of unity. \hbox{However},
$\mu_\alpha = - ((k_1^2 - 1)/({ k_1^1 +1})) - ({
2k_1}/({k_1^2+1})) i = q_1 + q_2 i$, with $q_1$ and $q_2$\break
\hbox{rational}. But only the first, second and fourth roots of
unity belong to $\Q(i)$, and $\mu_\alpha$ is none of these.
Finally, the case $n=1$ may be considered separately and shown not
to satisfy the boundary condition along $C_2$. Thus all $c_n=0$
and $\xi \equiv 0$. This completes the proof that no non-planar
holomorphic disks contribute to the Fukaya product of Lagrangian
planes.

\section{Relations}
\label{relsec}

In the previous section, we found that only planar disks
contribute to the product (\ref{m}). In this section, we use this
information to derive an explicit expression for the product and
some results about the structure of the ring~$\RR$.

The basic holomorphic disk with vertices $Y_0 \,{\in}\, L_0
\,{\cap}\, L_{k_1}$, $\rho^{k_1} (Y_1) \in \rho^{k_1}( L_0 \cap
L_{k_2})$ and $Y_2 \in L_0 \cap L_{k_1+k_2}$ has symplectic area
$({k_1 k_2}/({2(k_1+k_2)})) \omega (Y_1 - Y_0, f(Y_1 - Y_0)) $.
The other triangles relevant to the same composition are taken
into account by the generating function
\[
A_{Y_1 - Y_0}^{[k]} = \sum_{\lambda \in \Lambda_0} {\rm e}^{-\pi
k (\omega - i b)(Y_1 - Y_0 - \lambda, f(Y_1 - Y_0 - \lambda))}
\]
In the canonical lattice basis, these generating functions are
structure\break \hbox{constants} for the product
\begin{equation} \label{sympro}
Y_0 Y_1 := \sum_{Y_2 \in K(L_{k_2+k_1})_0}
A^{[k_1k_2/(k_2+k_1)]}_{Y_1 - Y_0} Y_2.
\end{equation}
This expression is compatible with Seidel's mirror map in the
following sense. Let $\widetilde T := (\Lambda_0 \oplus
\Lambda'^*)\otimes \R /(\Lambda_0 \oplus \Lambda'^*)$  be the
dual torus, mirror to  $T$. If $M$, $B$ and $N$ are as in
Section~\ref{sec1}, then $\widetilde T$ has complex structure
\[
\widetilde J:=
\left(
\begin{matrix}
1 & 0 \\[-1pt]
B & 1
\end{matrix}
\right)
\left(
\begin{matrix}
0 & -M^{-1}\\[-1pt]
M & 0
\end{matrix}
\right)
\left(
\begin{matrix}
1  &  0 \\[-1pt]
-B &  1
\end{matrix}
\right)
\]
while the mirror of the endomorphism $J_f$ is the polarization
\[
E_f := \left(
\begin{matrix}
0 & -N^{{\rm T}} \\[-1pt]
N & 0
\end{matrix}
\right),
\]
which is compatible with $\widetilde{J}$ by the symmetry of
$N^{{\rm T}} M$. Therefore, $E_f$ is the first  Chern class of
some line bundle (of zero characteristics) $\widetilde L$. Assume
that $\widetilde L$ is ample. We denote by $H_f$ the
Hermitian\vadjust{\pagebreak} form associated to $E_f$ and by
$S_f$ the $\C$-linear extension of $H_f$ restricted to
$(\Lambda'^*\otimes \R) \times (\Lambda'^*\otimes \R)$. Moreover,
let
\[
\widetilde J_b :=
\left(
\begin{matrix}
0 & -B^{-1} \\[-1pt]
B &  0
\end{matrix}
\right).
\]
Then, for all $x,y \in L_0$,
\[
(H_f-S_f)(x,y):=E_f(\widetilde J x, y)-iE_f(\widetilde J_b x, y)
= (\omega-i b)(x, f(y)) .
\]
This identification allows us to compare the Fukaya product with the classical product of theta functions.
We follow the treatment of \cite{B-L}. A canonical basis for $H^0(\widetilde T, \widetilde L)$ is given by the following theta functions:
\begin{align*}
\theta_c^{\widetilde L}(v) &:= {\rm
e}^{(({\pi}/{2})S_f(v,v)-({\pi}/{2})(H_f-S_f)(c+2v,c))}
\\[-2pt]
&\qquad\times \sum_{\lambda\in\Lambda(L_k)_0}{\rm e}^{\pi(H_f-S_f)
(c+v,\lambda) - ({\pi}/{2}) (H_f-S_f)(\lambda, \lambda)},
\end{align*}
where $v\in (\Lambda_0 \oplus \Lambda'^*)\otimes \R$ and $c\in K(L_0)_0$.  Moreover, the product
\[
\widetilde m{:}\ H^0(\widetilde T, \widetilde L) \otimes H^0
(\widetilde T,\widetilde L) \longrightarrow H^0(\widetilde T,
\widetilde L^2)
\]
can be expressed on this basis as
\begin{equation} \label{cpxpro}
\widetilde m (\theta_{c_1}^{\widetilde L} \otimes
\theta_{c_2}^{\widetilde L}) = \sum_{c_3\in K(L_2)_0}
\theta^{\widetilde L^2}_{c_2-c_1}(0) \theta^{\widetilde
L^2}_{c_3}
\end{equation}
Since the elements of the canonical bases for $\HOM _{D{\rm
Fuk}(T)}$ and $\HOM_{\mathcal D(\widetilde T)}$ are both labeled
by elements in $K(L)_0$, it is natural to put them in
correspondence $Y_c \leftrightarrow \theta_c$. In this way, we
identify the generators of the ring $\RR$ with the generators of
the homogenous coordinate ring of the projective embedding of the
mirror complex torus $\widetilde T$. Given a relation $ \sum
C^{c_1,\ldots, c_n} Y_{c_1} \cdots Y_{c_n} = 0$ in $\RR$, we can
reduce it to a linear one using the product~(\ref{sympro}).
Replacing $Y_c \leftrightarrow \theta_c$ everywhere yields another
valid linear relation, as both sets of generators are linearly
independent over $\C$. Finally, using the product (\ref{cpxpro})
and the identity
\[
\theta_c^{\widetilde L^{2k}}(0) = \sum_\lambda {{\rm e}}^{
-k\pi(H_f-S_f)(c-\lambda,c-\lambda)} = A^{[k]}_c,
\]
we can work backwards and obtain the relation $\sum
C^{c_1,\ldots, c_n} \theta_{c_1}, \ldots, \theta_{c_n} = 0$.
Clearly, this process can be reversed and we conclude that $\RR$
and the homogeneous coordinate ring of the embedded mirror are
isomorphic.

\enlargethispage{6pt}

At least in principle, we can apply the point of view of~[14] and present, using only the knowledge
of the Fukaya category of $T$, the mirror torus $\widetilde T$ as an explicit complete intersection
in some projective space, uniquely specified by the Lagrangian $L_0$ and the symplectomorphism
$\rho$. For example, a set of simple relations can~be recovered as follows. Let $m_k{:}\
\RR^{\otimes k} \to \RR$ be the $k$-fold multiplica\-tion. Then for any $z\in \RR$ of degree $k$
and for any pair $x_1\otimes\cdots\otimes x_k$, $y_1\,{\otimes}\,\cdots \,{\otimes}\, y_k$ such
that $x_i,y_i \in \RR$ all of degree 1 and\break $m_k(x_1\,{\otimes}\, \cdots \,{\otimes}\,
x_k)\kern-1pt=\kern-1pt\alpha_x z$, $m_k(y_1 \otimes \cdots \otimes y_k)\kern-1pt =\kern-1pt
\alpha_y z$, then $\alpha_y (x_1,\ldots, x_k)=\break\alpha_x(y_1,\ldots, y_k)$ is clearly a
\hbox{relation} in $\RR$. In the language of theta-functions, for $k=2,3$, these relations are
classically known as Riemann's theta relations and cubic theta relations, respectively. The linear
map $f$ can be chosen so that the above relations generate any other relation in $\RR$.

\enlargethispage{12pt}

\section{Dependence on $\rho$}

In this section, we investigate some examples and study how the
mirror family depends on the choice of symplectomorphism around
the large complex structure point.

\subsection{Twisted homogeneous coordinate rings}

So far, we assumed the symplectomorphism $\rho_f$ to be strictly
linear. \hbox{Relaxing} this condition and allowing affine
symplectomorphisms, one can reconstruct twisted homogeneous rings
as well. As an example, we now show that the noncommutative
projective plane can be reconstructed from a symplectic two-torus
together with the symplectomorphism $\rho(x,y) = (x+b, y+3x)$ for
any $b\neq \Z$. If $b \in \Z$ we recover the situation of~[14]
from which we adapt the notation as follows. Let $T = \R^2/\Z^2$
be the torus with coordinates $x, y$ and symplectic form $\tau dx
\wedge dy$. In the universal cover, we define three Lagrangians
$L_0 := \{y=0\}$, $L_1 := \{y = 3 (x-b)\}$ and $L_2 := \{y =
6(x-\break
 3b/2)\}$ for some $b \in \R$. Passing to the quotient, $L_0 \cap
L_1 = \{ X_i:=(i/3+b,0) | i = 0, 1, 2\} $, $L_0 \cap L_2 =
\{Y_j:=(j/6+3b/2,0) |  j = 0,\ldots , 5 \} $ and $L_1 \cap
L_2\break = \{ (k/3+2b, k+3b) |  k=0,1,2\} = \rho \{X_0, X_1,
X_2\}$. For general~\hbox{$b\in \R$}, the product formula of
\cite{Z} becomes
\[
X_i X_j = \sum_{k=0}^1 A_{i - j + 3k}(b) Y_{i + j + 3k} ,
\]
where we have put
\[
A_{k}(b) := \sum_n {{\rm e}}^{i\pi 6 \tau (n + k/6 + b/2)^2} =
\theta\left[\frac{k}{6}+\frac{b}{2},0\right](6\tau, 0).
\]
If $b \in \Z$, commutativity is ensured by the relation $A_k(b)= A_{6-k}(-b)$. For $b\notin \Z$, we get the relations
\begin{align*}
p X_2^2 + q X_0 X_1 + r X_1 X_0 & =  0 \\
p X_1^2 + q X_2 X_0 + r X_0 X_2 & =   0 \\
p X_0^2 + q X_1 X_2 + r X_2 X_1 & =   0
\end{align*}
where
\[
p := A_1 A_2 - A_4 A_5,\quad q:= A_3 A_4 - A_0 A_1, \quad r:= A_0
A_5 - A_3 A_2
\]
and $b$ dependence is understood. The noncommutative algebra,
which is the quotient of the noncommutative (associative)
homogeneous polynomial ring $\C\{X_0,X_1,X_2\}$ by the above
relations, is known as the Sklyanin algebra ${\rm Skl}_3(p,q,r)$.
If $b\in \Z$,  ${\rm Skl}_3(p,q,r)$ simply reduces to the ring of
homogeneous polynomials in three variables and so when $b\notin
\Z$, it makes sense to interpret it as the homogeneous coordinate
ring of a noncommutative projective plane. Numerical checks
confirm that elliptic curve $E$ with modular parameter $\tau$
(i.e., the mirror of $T$) has the equation
\[
Q_{pqr}:=X_0^3+X_1^3+X_2^3-\frac{p^3+q^3+r^3}{pqr}X_0X_1X_2=0.
\]
On the other hand, $Q_{pqr}$ generates the center of ${\rm
Skl}_3(p,q,r)$ and (as shown in \cite{ATV, BP, AKO})
\[
\frac{{\rm Skl}_3(p,q,r)}{Q_{pqr}} \cong \bigoplus_{n\ge 1}
H^0(E,L\otimes (t_b)^*L\otimes\cdots\otimes ((t_b)^{n-1})^*L)
\]
where $L$ is the degree 3 line bundle on $E$ which defines the
projective embedding such that $E={\rm
Proj}(C[X_0,X_1,X_2]/Q_{pqr})$ and $t_b\in {\rm Aut}(E)$ is the
translation by $b$. Moreover, the category of coherent sheaves on
$E$ embeds into the category of graded modules (up to torsion)
over ${\rm Skl}_3(p,q,r)$ which, in the language of \cite{ATV},
justifies the assertion that the mirror of $T$ sits as a
commutative curve into a noncommutative projective space.

\subsection{Quasihomogeneous coordinate rings} \label{qhom}

Consider the two-torus of \cite{Z} again, but with $\rho(x,y) =
(x,y+x)$ rather than the choice $(x,y+3x)$ which led to cubic
curves in the Hesse family.  The appearance of quasihomogeneous
coordinate rings will be quite natural.  To do the calculation,
we return to the philosophy of \cite{Z}, where the mirror map was
found without any prior knowledge of it.  The relations will be
shown to agree with Section~\ref{relsec}.

We have $L_k = \{y=kx\}.$  Let $X$ be the lattice point $0 \in
L_0\cap L_1$, $Y_0$\break and~$Y_1$~the points $Y_k = (k/2,0) \in
L_0\cap L_2$, and $Z_0, Z_1, Z_2$ the points\break $Z_k =
(k/3,0)\in L_0\cap L_3.$ Compute $X^2 = \aco{2}{0} Y_0 +
\aco{2}{1} Y_1$, where we define
\[
\aco{k}{j} = \theta\left[\frac{j}{k},0\right](k\tau,0), \quad
j\in \Z/k\Z.
\]
We put $Y = Y_0$ and $Z = Z_1.$  Then relations will necessarily
be quasihomogeneous with respect to the grading $|X| = 1$,
$|Y|=2$, $|Z|=3.$  Also note that using the (commutative) Fukaya
product, we can express $Z_k$ in terms of $X^3$, $XY$ and
$Z.$\footnote{Define
\[
\mu =
\frac{\aco{2}{0}\aco{6}{0}+\aco{2}{1}\aco{6}{3}}{\aco{2}{0}\aco{6}{2}
+\aco{2}{1}\aco{6}{1}},\quad \nu =
\frac{\aco{6}{0}}{\aco{6}{2}}.
\]
Then
\begin{gather*}
Z_0 = \frac{1}{(\aco{2}{0}\aco{6}{2}
+\aco{2}{1}\aco{6}{1})(\mu-\nu)}X^3 -
\frac{1}{\aco{6}{2}(\mu-\nu)}XY,\quad Z_1 = Z,
\\
\mbox{ and } Z_2 = -\frac{\nu}{(\aco{2}{0}\aco{6}{2}
+\aco{2}{1}\aco{6}{1})(\mu-\nu)}X^3 -
\frac{\mu}{\aco{6}{2}(\mu-\nu)}XY - Z. \end{gather*}} Similarly,
the results of \cite{Z} allow us to readily express the six
points $W_k = (k/6,0) \in L_0\cap L_6$ in terms of products
$Z_iZ_j.$  We only require $W_0 = a_1 Z_0^2 - a_3 Z_1 Z_2$ and
$W_2+W_4 = a_1(Z_1^2 + Z_2^2) - a_3(Z_0Z_1+Z_0Z_2)$,\break where
$a_k = \aco{6}{k}/(\aco{6}{0}\aco{6}{1}-\aco{6}{2}\aco{6}{3}).$
Writing the $W_k$ in terms of bilinears in $X^3$, $XY$ and $Z$
gives six of the seven quasihomogeneous monomials of degree six
in $X$, $Y$, and $Z$.  The remaining one is $Y^3$, which we
calculate in terms of the $W_k$ as
\[
Y^3 = \big(\aco{4}{0}\aco{12}{0}+\aco{4}{2} \aco{12}{6}\big)W_0 +
\big(\aco{4}{0}\aco{12}{4}+\aco{4}{2}\aco{12}{2}\big) (W_2+W_4).
\]
This gives a single relation in degree six.

Now given a polynominial equation as
\[
F\equiv Y^3-(p_0Z^2+p_1XYZ+p_2X^2Y^2+p_3X^3Z+p_4X^4Y+p_6X^6)=0,
\]
it is a simple matter to make ``linear'' changes of variables
(such as $Z \rightarrow\break Z + ({1}/{2p_0})(p_1XY + p_3 X^3)$)
to put the equation in a form with $p_0 = 1/4$, $p_1 = p_2 = p_3
= 0.$ Then working in an affine patch with coordinates $y =
Z/X^3$, $x = Y/X^2$, the equation has the Weierstrass form $y^2 =
4x^3 - g_2 x - g_3$,\footnote{Define $t_2  = p_2 - p_1^2/(4p_0)$,
$t_4 = p_4 - p_1p_3/(2p_0)$, $t_6 = p_6 - p_3^2/(4p_0).$ Then
$g_2^3 = (2/p_0)^2 \left[ t_4 + t_2^2/3 \right]^3$ and $g_3^2 =
(1/p_0)^2 \left[ t_6 + t_4 t_2/3 + 2t_2^3/27\right]^2.$} from
which one reads $j(\tau) = 1728 g_2^3 / (g_2^3 - 27g_3^2).$
Plugging in the power series as outlined above and in the
footnotes, we find (up to the first few dozen coefficients
checked by computer) the usual integer $q$-series expansion for
$j$, with $q = {\rm e}^{2\pi i \tau}$, i.e., the elliptic curve
has modular parameter $\tau$, and the usual mirror map is
established.

The ring $\mathcal R = \C[X,Y,Z]/F$ describes the mirror elliptic
curve $F=0$ inside weighted projective space $\P^2(1,2,3)$ (see
\cite{MR}). To compare with the\break \hbox{quadratic} relations
of Section~3, one considers the seven monomials
\[
V_0,\ldots,V_6 = X^6,X^4Y,X^3Z,X^2Y^2,XYZ,Y^3,Z^2
\]
as homogeneous projective coordinates on $\P^6$.  The nine
quadratic relations such as $V_0 V_4 = V_1 V_2$ describe the
(image of the) Veronese embedding $\P^2(1,2,3)\hookrightarrow
\P^6.$  In this description, the relation $F=0$ is linear in the
$V_k$, i.e., the projective line bundle has a unique section ---
the mirror of $X\in L_0\cap L_1.$

\subsection{Kummer varieties} \label{ku}

Seidel's method applies to quotients of abelian varieties as well.
By definition, a Kummer surface is the quotient of an abelian
variety with respect to the involution which reverses the
orientation of the lattice.  It is a singular surface with 16
singularities, and if it can be embedded as a hypersurface in
$\P^3$ it has equation~\cite{GD}
\begin{align*}
&A (X_0^4 + X_1^4 + X_2^4 + X_3^4)  + B (X_0^2 X_1^2 + X_2^2
X_3^2) + C (X_0^2 X_2^2 + X_1^2 X_3^2)
\\
&\quad + D (X_0^2 X_3^2 + X_1^2 X_2^2) + 2 E X_0 X_1 X_2 X_3 = 0
\end{align*}
We claim that we can reconstruct such a Kummer surface from a real
four-torus $T:=\R^4/\Z^4$ endowed with complex symplectic form
\[
\omega=\tau_1 dx_1 \wedge dy_1 + \tau_2 dx_2 \wedge dy_2 + \tau_3
(dx_1\wedge dy_2 + dx_2 \wedge dy_1),
\]
the standard involution $\iota(x_1, y_1, x_2, y_2) = - (x_1, y_1,
x_2, y_2)$ and the symplectomorphism $\rho(x_1, y_1, x_2, y_2) =
(x_1, y_1+ 2 x_1, x_2, y_2 + 2 x_2)$.  On $K:=T/\iota$,
\begin{gather*}
|L_0 \cap \rho L_0|=\binom{4}{1}\quad |L_0 \cap
\rho^2 L_0| = \binom{5}{2} \quad |L_0 \cap \rho^3 L_0| =
\binom{6}{3}
\\
|L_0\cap \rho^4 L_0| = \binom{7}{4} -1
\end{gather*}
so that, generically, among the homogeneous polynomials in $X_0, \ldots, X_3\in L_0 \cap \rho L_0$, we expect one relation in degree 4 to define the ring $\RR$ with no relations in degree 2 or 3. The computation of $\RR$ for general Kummer varieties will appear elsewhere \cite{A}.

To give an idea of the computation, here we illustrate the
degenerate case where $\tau:=\tau_1=\tau_2$ and $\tau_3 = 0$. In
particular, the matrices $M$, $B$ and $N$ (see
Section~\ref{sec1}) are simultaneously diagonalizable, a
situation that is mirror to the case of an abelian surface
polarized by the square of a reducible principal polarization.
Therefore we expect (e.g., from \cite{GD}, Proposition~4.23) a
map onto a quadric in $\P^3$. In the universal cover, consider
the Fukaya products
\[
\yco{k}{a,b}\yco{k}{c,d}:=\sum_{i,j\in \Z/2\Z\times \Z/2\Z}
\aco{2k}{c-a+2ki} \aco{2k}{d-b+2kj} \yco{2k}{a+c+2ki,b+d+2kj}
\]
where
\[
\yco{k}{a,b} : = \left(\frac{a}{2k},0,\frac{b}{2k},0\right)\in
L_0\cap \rho^k L_0
\]
In particular, defining
\[
X_0 := \yco{1}{0,0}, \quad X_1 := \yco{1}{1,0} , \quad X_2 :=
\yco{1}{0,1}  , \quad X_3 := \yco{1}{1,1},
\]
we have
\begin{equation} \label{ver}
X_0 X_3 =  \aco{4}{1} \aco{4}{1} \Big(\yco{2}{1,1} + \yco{2}{1,3}
+ \yco{2}{3,1} +\yco{2}{3,3}\Big)  = X_1 X_2
\end{equation}
The last equation defines the image of a Veronese embedding $\P^1
\times \P^1 \hookrightarrow \P^3$. On the other hand, $\rho$ can
be restricted to $\rho_i(x_i,y_i)=(x_i,y_i+2x_i)$, $i=1,2$. The
corresponding Seidel map $(\R^2/\Z^2,\omega_i = \tau dx_i \wedge
dy_i)\to \P^1$ can be used to find the equation of the mirror
elliptic curve in weighted projective space $\P(1,1,2)$, as in
the last section. The mirror map defined by $\rho$ is seen to be
the composition of the Veronese embedding and the Cartesian
product of the maps arising from $\rho_1$ and $\rho_2$.

\section*{Acknowledgments}

We would like to thank Paul Seidel for sharing his ideas. We are
grateful to The Fields Institute for partial support and for
hosting us while parts of this project were completed. The work
of E.Z. was supported in part by a Clay Senior Scholars
fellowship and by NSF grant DMS--0405859. Any opinions, findings
and conclusions or recommendations expressed in this material are
those of the authors and do not necessarily reflect the views of
the National Science Foundation (NSF).

\end{document}